\documentclass[11pt,reqno]{amsart}
 
 \usepackage{amsmath}
 \usepackage{amssymb}
\usepackage{bm}
 \usepackage{mathtools}
\usepackage{color}
 
\usepackage{graphicx}

\newcommand{\mysection}[1]{\section{#1}
      \setcounter{equation}{0}}

\newtheorem{theorem}{Theorem}[section]
\newtheorem{lemma}[theorem]{Lemma}

\newtheorem{corollary}[theorem]{Corollary} 

\theoremstyle{definition}

\theoremstyle{remark}
\newtheorem{remark}[theorem]{Remark}

\newcommand{\loc}{\text{\rm loc}}

 \makeatletter 
\makeatletter 
 \def\dashint{\operatorname{\,\,\,\mathclap{\int} \kern-.23em\text{\bf--}\!\!}}
 \makeatother

\def\dashnorm{\,\,\text{--}\kern-.4em\|}
\def\ninf{\qopname\relax\@empty{inf\phantom{p}\!\!\!}}

 \makeatother

\newcommand\bR{\mathbb{R}}

\newcommand\bC{\mathbb{C}}
\newcommand\bB{\mathbb{B}}

\begin{document}

\title[Two results related to elliptic and parabolic equations]
{Two results related to elliptic and parabolic equations. Case of Morrey spaces}

\author{N.V. Krylov}
 
\email{nkrylov@umn.edu}
\address{127 Vincent Hall, University of Minnesota,
 Minneapolis, MN, 55455, USA}


\dedicatory{To the memory of G. Da Prato}

\begin{abstract}
The paper presents a new short proof
of one of Adams's theorem and a $t$-trace-class theorem for parabolic Morrey
spaces.
 
\end{abstract}

\maketitle

\mysection{Introduction}

In this paper we discuss two results
from analysis which the author needed
in his investigation of strong solvability
of the It\^o stochastic equations with
singular coefficients. These results,
most certainly, are also useful in the theory
of partial differential equations.

The first result consists of a new
short proof of one of Adams's theorem,
which implies the possibility to estimate
the $L_{p}(\bR^{d})$-norm of $b^{i}D_{i}u$
through the $L_{p}(\bR^{d})$-norm of $D^{2}u$
when $p<d$ and $b\not\in L_{d}(\bR^{d})$,
so that Sobolev embedding theorems
do not guarantee that $b^{i}D_{i}u
 \in L_{p}(\bR^{d})$. The results of that kind first appeared in Fefferman
\cite{Fe_83} (1983) in PDE terms when $p=2$.
Then Adams \cite{Ad_86} (1986) proved them
in terms of Riesz potentials for any
$p>1$, absolutely differently
from \cite{Fe_83} and, probably, without
even knowing that it existed.
The Adams theorem allows us to estimate
$b^{i}D_{i}u$ by    fractional
powers   of the Laplacian of $u$.
After that Chiarenza-Frasca \cite{CF_90}
(1990)
wrote a paper generalizing the result
of Fefferman about estimates of $b^{i}D_{i}u$ to any $p>1$ using a clever PDE
argument and, probably, without
 knowledge of \cite{Ad_86}.

The author gave two poofs of the Adams
theorem. One generalizes it to
parabolic setting \cite{4}
  and the other \cite{8} (2024) gives a shorter
proof  of the   Adams theorem
based on the ideas from \cite{CF_90}. Still the proof in \cite{8} uses some rather
deep results from Real Analysis.
Here, following \cite{4}, we give
even shorter proof of the  Adams theorem
based on only well-know facts from Real Analysis such as the Hardy-Littlewood
maximal function inequality and the
Fefferman-Stein sharp function theorem.
These, one can find, for instance, in
\cite{Kr_08}.

The second result of this paper is
about the $t$-traces of parabolic 
Morrey spaces of functions $u(t,x)$
which have one derivative in $t$
and two derivatives in $x$.
This result is aimed at implying the space differentiability at each time $t$
of the evolutional family related
to parabolic PDEs with singular first-order coefficients when the Sobolev-space
approach does not work.
As far as the author is aware, there are no previous results on this subject--one more evidence that the Morrey space
theory of PDEs is underdeveloped.
What we obtain is, certainly, not
the sharpest result (the sharpest one should, probably,
involve some kind of Besov-Morrey space)
but it is sufficient for the applications
we have in mind.

\mysection{Adams's theorem}

Here we present a proof of Theorem 1 of \cite{8}, which is a particular case of Theorem 7.3 of 
\cite{Ad_86} or, if $\alpha=1$, of  the Theorem
in \cite{CF_90} (where $d\geq3$).  The need in such results arises when one wants to
add to $\Delta$ the first order terms
$b^{i}D_{i}$ and wants to solve the corresponding elliptic equation in
$W^{2}_{p}$ with $p<d$ and $b\not\in L_{d}$ but having lower integrability.
It turns out that, if
$b$ belongs to the Morrey space $\dot E_{q,1}$ with $p<q$, then
$b^{i}D_{i}$ is a bounded operator
from $W^{2}_{p}$ to $L_{p}$ (see Corollary \ref{corollary 3.12.1}).

For $d\geq2$ let $\bR^{d}$ be the Euclidean space of points $x=(x^{1},...,x^{d})$. 
Define $B_{\rho}(x)=\{y:|y-x|<\rho\}, 
B_{\rho}=B_{\rho}(0)$, and let $\bB_{\rho}$ be the collection of $B_{\rho}(x)$.
Denote $D_{i} =\partial /\partial x^{i}$,
$Du=(D_{i}u), D_{ij}=D_{i}D_{j}$, $D^{2}u
=(D_{ij}u)$. 
For measurable $\Gamma\subset \bR^{d}$ set $|\Gamma|$ to
be the volume of $\Gamma$ and introduce
$$
\dashint _{\Gamma}f\,dx=\frac{1}{|\Gamma|}
\int _{\Gamma}f\,dx,\quad 
\dashnorm f\|^{p}_{L_{p}(\Gamma)}
=\dashint_{\Gamma} |f|^{p}\,dx,
$$
$$
M_{\beta}g(x)=\sup_{\rho>0}\rho^{\beta}\dashnorm g\|_{L_{1}(B_{\rho}(x))},\quad
u^{\#}(x)=\sup_{\rho>0}\dashint_{B_{\rho}(x)}
\dashint_{B_{\rho}(x)}|u(y)-u(z)|\,dydz.
$$

Let
$\beta\geq 0,q\in[1,\infty)$. In our system of notation the Morrey space $\dot E _{q,\beta}$
is the collection of functions such that
$$
\|f\|_{\dot E_{q,\beta}}:=\sup_{\rho>0} \rho^{\beta}
\sup_{ B\in\bB_{\rho} }
\dashnorm f\|_{L_{q}(B)}<\infty.
$$
Also for $0<\alpha<d$ set
$$
P_{\alpha}f(x)=\int_{\bR^{d}}\frac{f(y)}{|x-y|^{d-\alpha}}\,dy,
$$
which is the Riesz potential of $f$.

\begin{theorem}
                \label{theorem 3.11.1}
Let $1<p<q$,   $b\in \dot E_{q,\alpha},f\in L_{p}$. Then
\begin{equation}
                            \label{3.12.2}
\|bP_{\alpha}f\|_{L_{p}}\leq N(d,\alpha,p,q)\|b\|_{ \dot E_{q,\alpha}}\|f\|_{L_{p}}.
\end{equation}
\end{theorem}

To prove the theorem we need
the following result, taken from
\cite{Ad_75} where the proof is slightly  different.
\begin{lemma}
                     \label{lemma 3.11.2}
There is $N=N(d,\alpha)$ such that for any bounded $g$ with compact support
\begin{equation}
                            \label{3.12.1}
P^{\#}_{\alpha}g:=(P_{\alpha}g)^{\#}\leq N M_{\alpha}g.
\end{equation}
 
\end{lemma}

Proof. Introduce
$g=f+h$, $f=gI_{B_{2}}$. Observe that for $|x|\leq 1$ we have
$$
|DP_{\alpha}h(x)|\leq N\int_{B^{c}_{2}}\frac{g(y)}
{|y-x|^{d-\alpha+1}}\,dy\leq N
\int_{B^{c}_{2}}\frac{g(y)}
{|y|^{d-\alpha+1}}\,dy
$$
$$
=N\int_{2}^{\infty}\frac{1}{r^{d-\alpha+2}}
\int_{B_{r}\setminus B _{2}}g(y)\,dy\,dr\leq
N\int_{2}^{\infty}\frac{1}{r^{d-\alpha+2}}
\int_{B_{r} }g(y)\,dy\,dr
$$
$$
\leq NM_{\alpha}g(0)\int_{2}^{\infty}r^{-2}\,dr=NM_{\alpha}g(0).
$$

Also
$$
\int_{B_{1}}\int_{B_{1}}|P_{\alpha}f(x)-
P_{\alpha}f(y)|\,dxdy\leq N
\int_{B_{1}}P_{\alpha}f(x)\,dx
$$
$$
=N\int_{\bR^{d}}\frac{1}{|y|^{d-\alpha}} \Big(\int_{B_{1}}
 f(x-y) \,dx\Big)\,dy=N\int_{B_{3}}\frac{1}{|y|^{d-\alpha}} \Big(\int_{B_{1}}
 f(x-y) \,dx\Big)\,dy
$$
$$
\leq N\int_{B_{3}}\frac{1}{|y|^{d-\alpha}}
\Big(\int_{B_{4}}
 f(x ) \,dx\Big)\,dy\leq NM_{\alpha}g(0).
$$

It follows that
$$
\dashint_{B_{1}}\dashint_{B_{1}}|P_{\alpha}g(x)-
P_{\alpha}g(y)|\,dxdy \leq N M_{\alpha}g(0)
$$
and by self-similarity we get \eqref{3.12.1}
at the origin. In any other point the proof
is similar. The lemma is proved. \qed

{\bf Proof of Theorem \ref{theorem 3.11.1}}. 
We may assume that  $b,f\geq0$ and let them
be bounded with compact support. Set $u=P_{\alpha}f$ and observe that
$$
I:=\int_{\bR^{d}} b^{p}u ^{p}\,dx=
\int_{\bR^{d}}\Big(P_{\alpha}( b^{p}u ^{p-1})\Big)f\,dx\leq \|f\|_{L_{p}}
\|P_{\alpha}( b^{p}u ^{p-1})\|_{L_{p'}}
$$
$$
\leq N\|f\|_{L_{p}}
\|P^{\#}_{\alpha}( b^{p}u ^{p-1})\|_{L_{p'}},
$$
where $p'=p/(p-1)$ and the last inequality follows by
the Fefferman-Stein theorem (see \cite{FS_82}
  or Theorem 3.2.10 in \cite{Kr_08}).
By Lemma \ref{lemma 3.11.2} we can replace
$P^{\#}_{\alpha}$ with $M_{\alpha}$. After that note that
$$
\rho^{\alpha}\dashint_{B_{\rho}(x)}
b^{p}u ^{p-1}\,dy=
\rho^{\alpha}\dashint_{B_{\rho}(x)}
b (bu) ^{p-1}\,dy\leq \|b\|_{E_{q,\alpha}}
M^{1/q'}_{0}[(bu)^{(p-1)q'}](x),
$$
which implies that
$$
M_{\alpha}( b^{p}u ^{p-1})\leq \|b\|_{E_{q,\alpha}}
M^{1/q'}_{0}[(bu)^{(p-1)q'}].
$$
Since $p'/q'>1$, by the Hardy-Littlewood
maximal function theorem
$$
\|M^{1/q'}_{0}[(bu)^{(p-1)q'}]\|_{L_{p'}}
\leq N\| bu\|_{L_{p}}^{p-1}=NI^{1-1/p}.
$$

It follows that $I\leq N\|b\|_{E_{q,\alpha}}\|f\|_{L_{p}}I^{1-1/p}$, which yields \eqref{3.12.2}
since $I<\infty$ in light of the fact that  $b$ and $f$ are bounded and have compact support. After that
approximations and
Fatou's lemma take  care of the general case.
The theorem is proved. \qed

\begin{corollary}
                  \label{corollary 3.12.1}
If $u\in C^{\infty}_{0}$ and $1<p<q $, then
$$
\|bu\|_{L_{p}}\leq N(d,p,q)\|b\|_{\dot E_{q,1}}
\|Du\|_{L_{p}},\quad \|bDu\|_{L_{p}}\leq N(d,p,q)\|b\|_{\dot E_{q,1}}
\|D^{2}u\|_{L_{p}}
$$
\end{corollary}

Indeed, as is well known,
$|u|\leq NP_{1}|Du|$.

\begin{corollary}[Hardy's inequalities]
                \label{corollary 11.28.1}
If $u\in C^{\infty}_{0}$ and $1<p<d $, 
$1<r<d/2$, then
$$
\int_{\bR^{d}}\frac{|u|^{p}}{|x|^{p}}
\,dx\leq N\|Du\|_{L_{p}}^{p},
\quad \int_{\bR^{d}}\frac{|u|^{r}}{|x|^{2r}}
\,dx\leq N\|D^{2}u\|_{L_{r}}^{r}.
$$
\end{corollary}

Indeed, as is easy to see, $|x|^{-1}$
is of class $E_{q,1}$ and $|x|^{-2}$
is of class $E_{q/2,2}$ for any
$q\in(1,d)$.

\begin{remark}
                   \label{remark 3.12.1}
The statements of Theorem \ref{theorem 3.11.1}
and Corollary \ref{corollary 3.12.1}
are meaningful only if $q\leq d$ because
otherwise $\|b\|_{E_{q,\alpha}}<\infty$ only if $b=0$.
\end{remark}

From the point of view of the theory
of elliptic equations the requirement
$b\in \dot E_{q,1}$ seems too restrictive
because it excludes, for instance,
bounded $b$. Therefore, the following
condition looks reasonable
$$
\hat b_{p_{b},\rho_{b}}:=
\sup_{\rho\leq \rho_{b}}\rho
\sup_{B\in \bB_{\rho}}\dashnorm b\|_{L_{p_{b}}(B)}<\infty,
$$
where $p_{b}\in (1,\infty)$ and $\rho_{b}
\in (0,\infty)$ are some numbers.

\begin{lemma}
                \label{lemma 7,10.3}
(i) If $p_{b}\leq d$, then for any
$r>0$ and
$B\in \bB_{r}$ we have $
r\dashnorm I_{B_{\rho_{b}}}b\|_{L_{p_{b}}(B)}\leq \hat b_{p_{b},\rho_{b}}$.

(ii) If $p_{b}\geq d$, then 
$$
\hat b_{p_{b},\rho_{b}}=\rho_{b}\sup_{B\in \bB_{\rho_{b}}}
\dashnorm b\|_{L_{p_{b}}(B)}=
N(d,p)\rho^{1-d/p_{b}}_{b}\sup_{B\in \bB_{\rho_{b}}}
\| b\|_{L_{p_{b}}(B)}.
$$ 
\end{lemma}

Proof. (i) If $r\leq \rho_{b}$,
then the assertion follows from the definition of $\hat b_{p_{b},\rho_{b}}$.
However, if $r>\rho_{b}$, then
$$
r\dashnorm I_{B_{\rho_{b}}}b\|_{L_{p_{b}}(B)}=r^{1-d/p_{b}}|B_{1}|^{-d/p_{b}}
\|I_{B_{\rho_{b}}\cap B}b\|_{L_{p_{b}}}
$$
$$
\leq \rho_{b}^{1-d/p_{b}}|B_{1}|^{-d/p_{b}}
\|I_{B_{\rho_{b}} }b\|_{L_{p_{b}}}
=\rho_{b}\dashnorm b\|_{L_{p_{b}}(B_{\rho_{b})}}\leq \hat b_{p_{b},\rho_{b}}.
$$

Assertion (ii) is obvious and the lemma is proved. \qed

\begin{remark}
             \label{remark 7,12.1}
The common requirement is for $\hat b_{p_{b},\rho_{b}}$ to be ``small enough''. Therefore, it is useful to know that,
if $p_{b}>d$ and $\sup_{B\in \bB_{1}}
\|b\|_{L_{p_{b}}(B)}<\infty$, then
$\hat b_{p_{b},\rho_{b}}$ can be made as small as we wish on account of
taking $\rho_{b}$ small enough.
\end{remark}

The following
 generalizes
and implies Corollary  \ref{corollary 3.12.1}
as $\rho_{b}\to\infty$.

\begin{lemma}   
                         \label{lemma 11.4.100}
For $1<p<p_{b} $ there exists a constant $N=N(p,p_{b},  d )$ such that
\begin{equation}
                       \label{11.4.30}
\int_{\bR^{d}}|b(x)|^{p}|u(x)|^{p}\,dx
\leq N \hat b_{p_{b},\rho_{b}}^{p}
\int_{\bR^{d}} |D u|^{p}\,dx
+N \rho_{b}^{-p }\hat b_{p_{b},\rho_{b}}^{p} 
\int_{\bR^{d}} | u|^{p}\,dx
\end{equation}
as long as $u\in C^{\infty}_{0}$.
\end{lemma}

 Proof.  First, observe that,
as it is not hard to see, 
\eqref{11.4.30} is scale invariant.
Therefore, we may assume that $\rho_{b}=1$.

{\em Case $p_{b}\leq d$\/}.
Take $\zeta\in C^{\infty}_{0}(B_{1})$,
$\zeta\geq0$, such that
\begin{equation}
                         \label{11.4.5}
\int_{B_{1}}\zeta^{ p}\,dx=1, \quad
  |D\zeta|\leq N(d) I_{B_{1}},
\end{equation}
and set $\zeta_{y}(x)=\zeta(x+y)$.
Then owing to Corollary  \ref{corollary 3.12.1} and 
Lemma \ref{lemma 7,10.3}, for any
$y\in\bR^{d}$
$$
\int_{\bR^{d}}|b |^{p}|\zeta_{y} u|^{p}\,dx=\int_{\bR^{d}}|bI_{B_{1}(y)} |^{p}|\zeta_{y} u|^{p}\,dx
\leq N \hat b_{1}^{p}
\int_{\bR^{d}}|D(\zeta_{y}u)|^{p}\,dx
$$
$$
=N \hat b_{1}^{p}
\int_{\bR^{d}}\zeta^{p}_{y}|Du|^{p}\,dx+N \hat b_{1}^{p}
\int_{\bR^{d}} |D
\zeta_{y}|^{p}\,|u|^{p}\,dx.
$$
By integrating through these
inequalities over $y\in\bR^{d}$
and using \eqref{11.4.5} we come to 
\eqref{11.4.30}.

{\em Case $p_{b}> d$\/}. By H\"older's
inequality and embedding theorems,
for $r=pp_{\beta}/(p_{\beta}-p)$
$$
\int_{\bR^{d}}|b |^{p}|\zeta_{y} u|^{p}\,dx\leq\|b\|^{p}_{L_{p_{b}(B_{1}(y))}}\|\zeta_{y} u\|^{p}
_{L_{r}}\leq N\hat b_{1}^{p} (\|D(\zeta_{y} u)\|^{p }_{L_{p}}
+\| \zeta_{y} u \|^{p }_{L_{p}})
$$
$$
\leq N \hat b_{1}^{p}
\int_{\bR^{d}}\zeta^{p}_{y}|Du|^{p}\,dx+N \hat b_{1}^{p}
\int_{\bR^{d}} |D
\zeta_{y}|^{p}\,|u|^{p}\,dx,
$$
and we finish the proof as in the first case. The lemma is proved. \qed

\begin{remark}
                 \label{remark 11.28.1}
From the point of view of the theory of elliptic equations
the most desirable version of \eqref{11.4.30}
would be, of course, 
\begin{equation}
                            \label{10,26.2}
\|b|Du|\,\|_{L_{p}}
\leq \varepsilon\| D^{2}u\|_{L_{p}}+N(\varepsilon)
\|u\|_{L_{p}}
\end{equation}
for any $\varepsilon>0$ with   $N(\varepsilon)$ independent
of $u$. It turns out that
this is impossible for $\varepsilon$ small if $p<d$ and $\hat b_{p_{b},\rho_{b}}<\infty$
for $p<p_{b}<d$.

To show that let $p\in(1,d )$  and $f(t)$, $t\geq0$, be a smooth function,
equal to zero for $t\leq 1$, equal
to $t $ for $t\in [2,3]$
and equal to zero for $t\geq4$. Clearly,
$|f'(t)|,|f''(t)|\leq N $, for some
constants $N$ and $f'(t)=1$ for $t\in [2,3]$. Now, for small $\kappa\in(0,1/4 ]$ define $u_{\kappa}(x)=f(|x|/\kappa)$.
Then
$$
D_{i}u_{\kappa}=f'\frac{x^{i}}{
\kappa|x|},\quad D_{ij}u_{\kappa}=
f''\frac{x^{i}x^{j}}{\kappa^{2}|x|^{2}}+f'\frac{1}{\kappa|x|}
\Big(\delta^{ij}-\frac{x^{i}x^{j}}{
|x|^{2}}\Big).
$$
It follows that $|D^{2}u_{\kappa}
|\leq N\kappa^{ -2} $
in $B_{1}$ and $|Du|=\kappa^{-1}$
in $B_{3\kappa}\setminus B_{2\kappa}$.
With $b=|x|^{-1}$ this yields
$$
\dashnorm bDu_{\kappa}\|_{L_{p} }
\geq N_{1} \kappa ^{-2},\quad
\dashnorm D^{2}u_{\kappa}\|_{L_{p} }
\leq N_{2}\kappa^{ -2}  .
$$
By adding to this that $\dashnorm u_{\kappa}\|_{L_{p} }
\leq N\sup |f|$, we see that \eqref{10,26.2}
can only hold if $N_{1} \leq\varepsilon N_{2}$.

\end{remark}

\mysection
{Trace theorem for parabolic Morrey spaces}

                \label{section 3.27.1}

Define $C_{\rho}=[0,\rho^{2})\times B_{\rho}$ and for $(t,x)\in \bR\times\bR^{d}=\bR^{d+1}$ let $C_{\rho}(t,x)=(t,x)+C_{\rho}$ and $\bC_{\rho}$ be the collection
of $C_{\rho}(t,x)$.
Take $p,q\in(1,\infty)$ and
   $\beta\geq 0$ and introduce the
parabolic
Morrey space $E_{p,q,\beta} $
as the set of $g\in  L_{p,q,\loc}$ 
\index{$A$@Sets of functions!$E_{p,q,\beta}$}%
such that  
\begin{equation}
                             \label{8.11.02}
\|g\|_{E_{p,q,\beta} }:=
\sup_{\rho\leq 1,C\in\bC_{\rho}}\rho^{\beta}
\dashnorm g  \|_{ L_{p,q}(C)} <\infty ,
\end{equation}  
where
$$
\dashnorm  \cdot \|_{ L_{p,q}(\Gamma)}:=
\|I_{\Gamma}\|_{ L_{p,q}}^{-1}\|\cdot\|_{ L_{p,q}(\Gamma)},\quad \|g\|_{ L_{p,q}(\Gamma)}^{q}=
\int_{\bR}\Big(\int_{\bR^{d}}|g|^{p}\,dx
\Big)^{q/p}\,dt.
$$
Define $\partial_{t}=\partial/\partial t$ and
$$
E^{1,2}_{p,q,\beta} =\{u:u,Du,D^{2}u,
\partial_{t}u\in E_{p,q,\beta} \},
$$
where    $Du,D^{2}u,
\partial_{t}u$ are Sobolev derivatives,
and 
\index{$A$@Sets of functions!$E^{1,2}_{p,q,\beta}$}%
provide $E^{1,2}_{p,q,\beta} $ with an obvious norm.
The subsets of these spaces
consisting of functions independent of
$t$ is denoted by $E_{p,\beta}$
and $E^{2}_{p,\beta}$, respectively.

Fix a nonnegative $\zeta\in C^{\infty}_{0}
(\bR^{d+1})$ with unit integral and
for $\varepsilon>0$ set
$$
\zeta_{\varepsilon}(t,x)=\varepsilon^{-d-2}
\zeta(t/\varepsilon^{2},x/\varepsilon),
\quad u^{(\varepsilon)}=u*\zeta_{\varepsilon}.
$$

In Real Analysis traditionally
instead of $E_{p,q,\beta}$ one considers its ``homogeneous'' version
$\dot E_{p,q,\beta}$ 
defined 
\index{$A$@Sets of functions!$\dot E_{p,q,\beta}$}%
as the set of $g$ such that
$$
\|g\|_{\dot E_{p,q,\beta} }:=
\sup_{\rho<\infty,C\in\bC_{\rho}}\rho^{\beta}
\dashnorm g  \|_{ L_{p,q}(C)} <\infty 
$$
(note $\rho<\infty$).

Our goal in this section is to show that the $t$-traces
of functions in $E^{1,2}_{p,q,\beta}$
possess some regularity as $L_{p}$-functions. For $\gamma=0$ or $1$ set
$$
D^{\gamma}=D\quad\text{if}\quad
\gamma=1\quad\text{and}\quad D^{\gamma}=1\quad
\text{if}\quad \gamma=0.
$$
Below by
$D^{\gamma}u(0,\cdot)$ we mean the limit
in $L_{r}(B)$ for any ball $B$
of $D^{\gamma}u^{(\varepsilon)}(0,\cdot)$
as $\varepsilon\downarrow 0$.
The existence of this limit easily follows from Lemma \ref{lemma 6,17.1} and Corollary \ref{corollary 6,19.1} below.

\begin{theorem}
              \label{theorem 6,6,1}
Take  
$r\in[ p,\infty)$, $\mu>0$    and assume
that  
$$
2-\gamma< \beta \leq \frac{d}{p}+\frac{2}{q}
<2-\gamma+\frac{d}{r},
\quad \kappa:=\gamma+\frac{d}{p}
+\frac{2}{q}-\frac{d}{r}\leq \mu
<2 .
$$

 Then for any
$u\in E_{p,q,\beta}$ the trace $D^{\gamma}u(0,\cdot)$ is uniquely defined
and for any $\varepsilon>0$
\begin{equation}
                     \label{6,6.1}
\|D^{\gamma}u(0,\cdot)\|_{E_{r,\beta+\gamma-\mu}(\bR^{d})}\leq N\varepsilon \|\partial_{t}u,
D^{2}u\|
_{ E _{p,q,\beta}}
+N\varepsilon^{-\mu/(2-\mu)}
\|u\|
_{ E _{p,q,\beta}},
\end{equation} 
\begin{equation}
                     \label{8,27.1}
\|D^{\gamma}u(0,\cdot)\|_{E_{r,\beta+\gamma-2}(\bR^{d})}\leq N  \|u\|
_{ E^{1,2} _{p,q,\beta}},
\end{equation} 
where  the constants $N$ depend only
on $d,p,q,\beta,\mu$.

\end{theorem}

\begin{remark}
               \label{remark 11.29.1}
Obviously, $E^{1,2}_{p,q,\beta}\subset
W^{1,2}_{p,q,\loc}$ and one can show that
$E^{1,2}_{p,q,\beta}\not\subset
W^{1,2}_{p+\varepsilon,q+\varepsilon,\loc}$ no matter how small $\varepsilon>0$ is.
Therefore, in terms of local summability of derivatives, general functions in 
$E^{1,2}_{p,q,\beta}$ are not much better
than $W^{1,2}_{p,q,\loc}$. For the latter
class the trace theorems for $\gamma=1$ (see Lemma 
\ref{lemma 6,17.1}) , basically, only can guarantee that $Du(0,\cdot)\in L_{r,\loc}$
and, if $r<d$, this does not yield even the boundedness of $u(0,\cdot)$. At the same time \eqref{8,27.1} and the Morrey
theorem (see, for instance, Theorem 10.2.1
in \cite{Kr_08}) imply that $u(0,\cdot)$
is $2-\beta$ H\"older continuous, provided
$\beta<2$, that is, almost Lipschitz continuous. Of course, this is at the expense of $u\in E^{1,2}_{p,q,\beta}$.

\end{remark}

\begin{remark}
               \label{remark 8,25.1}
From the probability   point of
view  the most important particular case
of \eqref{6,6.1} is when $\gamma=1$,
$r=p$ (so that $q>2$), $\mu=\kappa$,
and 
$$
\frac{d}{p}+\frac{2}{q}\geq \beta>1 .
$$
 
 In that case
for any
$u\in E^{1,2}_{p,q,\beta}$  and any $\varepsilon>0$ (observe that $\beta>\beta+1-\kappa$)
\begin{equation}
                     \label{6,6.10}
\|D u(0,\cdot)\|_{E_{p,\beta } }\leq N\varepsilon \|\partial_{t}u,
D^{2}u\|
_{ E _{p,q,\beta}}
+N\varepsilon^{-(q+2)/(q-2)}
\|u\|
_{ E _{p,q,\beta}}.
\end{equation} 
 
\end{remark}

To prove Theorem \ref{theorem 6,6,1},
first, we need the following corollary of
Theorem 10.2 of \cite{BIN_75} which
we give with a different proof for completeness.

\begin{lemma}
                 \label{lemma 6,17.1}
Let  
$r\geq p$  and assume
that $\kappa<2$.
 Then for any
$u\in W^{1,2}_{p,q}$ and $\varepsilon>0$ we 
have
\begin{equation}
                            \label{6,17.2}
\|D^{\gamma}u(0,\cdot)\|_{L_{r}}\leq
N\varepsilon \|\partial_{t}u,D^{2}u
\|_{L_{p,q} }+N\varepsilon^{-\kappa/(2-\kappa)}
\|u
\|_{L_{p,q} }.
\end{equation}
\end{lemma}

Proof. The case of arbitrary $\varepsilon>0$
is reduced to that of $\varepsilon=1$
by using   self-similarity.
To treat $\varepsilon=1$ take $\zeta\in
C^{\infty}_{0}(\bR)$ such that $\zeta(t)=1$
for $|t|\in[0,1]$, $\zeta(t)=0$ for $|t|\geq 2$, and define $-f=
\partial_{t}(\zeta u)+\Delta(\zeta u)$. Introduce
$$
P_{\gamma}(t,x)=t^{-(d+\gamma)/2}e^{-|x|^{2}/(8t)}.
$$
We know   that $\zeta u=Rf$, where
$$
Rf(t,x)=N(d)\int_{0}^{\infty}s^{-d/2}
\int_{\bR^{d}}e^{-|x-y|^{2}/(4s)}
f(t+s,y)\,dyds.
$$
It follows that ($y^{0}:=1,y^{1}:=y$)
$$
D^{\gamma}u(0,x)= -N(d)\int_{0}^{\infty}t^{-d/2}
\int_{\bR^{d}}\Big(\frac{y}{2t}\Big)^{\gamma}e^{-| y|^{2}/(4t)}f(t,x-y)\,dydt.
$$
By observing that $|y/\sqrt t|
e^{-| y|^{2}/(4t)}\leq Ne^{-| y|^{2}/(8t)}$ we conclude that
$$
|D^{\gamma}u(0,x)|\leq NF(x),
$$
where
$$
 F(x) =
\int_{0}^{2} 
\int_{\bR^{d}}P_{\gamma}(t,y)|f(t,x-y)|\,dydt .
$$

By Minkowski's inequality
$$
\|F\|_{L_{r} }
\leq \int_{0}^{2}\Big(
\int_{\bR^{d}}\Big(\int_{\bR^{d}}
P_{\gamma}(t,x-y)|f(t,y)|\,dy\Big)^{r}
\,dx\Big)^{1/r}
\,dt,
$$ 
where inside the integral with respect
to $t$ we have the norm of   convolution, so that by Young's
inequality this expression is dominated by
$$
\|f(t,\cdot)\|_{L_{p}}\|P_{\gamma}(t,\cdot)
\|_{L_{s}},
$$
where $1/s=1+1/r-1/p$ ($ \leq 1$ since $r\geq p$). An easy computation
shows that $\|P_{\gamma}(t,\cdot)
\|_{L_{s}}=N(d)t^{\alpha}$ with
$\alpha=-\gamma/2+(d/2)(1/r-1/p)$, which  yields
$$
\|F\|_{L_{r} }
\leq N\int_{0}^{2}\|f(t,\cdot)\|_{L_{p} }t^{\alpha}\,dt.
$$
Now use H\"older's inequality along with the observation that
$\alpha q/(q-1)>-1$, due to 
the assumption that $\kappa<2$, that is,  $2-\gamma +d/r>
d/p+2/q$, to conclude that
$$
\|F\|_{L_{r} }
\leq N \|f \|_{L_{p,q}  }.
$$ 
The lemma is proved.\qed

\begin{corollary}
           \label{corollary 6,19.1}
For any $\rho\leq 1,\varepsilon>0$ and
$u\in W^{1,2}_{p,q}(C_{2\rho})$ we have
$$
\dashnorm D^{\gamma} u(0,\cdot)\|_{L_{r}(B_{\rho})}
\leq N\varepsilon  
\dashnorm \partial_{t}u,D^{2}u
\|_{L_{p,q}(C_{2\rho})}
$$
$$
+N
(\varepsilon \rho^{-2}+\varepsilon^{-\kappa/(2-\kappa)}\rho^{(2\kappa-2\gamma)/(2-\kappa)})
\dashnorm u
\|_{L_{p,q}(C_{2\rho})}.
$$
\end{corollary}

Indeed, the case of $\rho<1$ is reduced
to $\rho=1$ by means of parabolic dilation. In the latter case
take $\zeta\in C^{\infty}_{0}
(\bR^{d+1})$ such that $\zeta=1 $ on $C_{1}$ and $\zeta=0$ in $\bR^{d+1}_{0}
\cap C_{2}$. Then  use
\eqref{6,17.2} to see that
$$
\|D^{\gamma}u(0,\cdot)\|_{L_{r}(B_{1})}
\leq N\varepsilon  
\|\partial_{t}(\zeta u),D^{2}(\zeta u)
\|_{L_{p,q} }+N\varepsilon^{-\kappa/(2-\kappa)}
\|u
\|_{L_{p,q}(C_{2}) }
$$
$$
\leq N\varepsilon  
\|\partial_{t}u,D^{2}u
\|_{L_{p,q}(C_{2}) }+N\varepsilon \|u,Du\|_{L_{p,q}(C_{2}) }
+N\varepsilon^{-\kappa/(2-\kappa)}
\|u
\|_{L_{p,q}(C_{2}) }.
$$
After that it only remains to use
the classical Sobolev-space interpolation inequality
$$
\|Du(t,\cdot)\|_{L_{p}(B_{2})}\leq
N \|D^{2}u(t,\cdot)\|_{L_{p}(B_{2})}+
N\| u(t,\cdot)\|_{L_{p}(B_{2})}.
$$
 
In the following lemma $\gamma$ can be any
number in $[0,d+2)$.
\begin{lemma}
                  \label{lemma 11.28.1}

Let $0\leq\gamma<\beta\leq d+2$. Then there exist   constants $N$
($<\infty$)  such that for any $f\geq0$
and $\rho\in(0,\infty)$ we have
\begin{equation}
                             \label{1.17.2}
P_{\gamma}(I_{C^{c} _{ \rho}}f)(0)
\leq N\rho^{\gamma-\beta}M_{\beta} 
f(0) .
\end{equation}
\end{lemma}

Proof. Clearly, we may assume that $f$ is bounded
with compact support.  Set $Q_{1}=\{(s,y):|y|\geq\sqrt s\}$,
$Q_{2}=\{(s,y):|y|\leq \sqrt s\}$.
Dealing with $P_{\gamma}(fI_{Q_{1}})$
we observe that $p_{\gamma}(s,r)\leq Nr^{-(d+2-\gamma)}$ if $r\geq\sqrt s$. Therefore,
$$
P_{\gamma}(fI_{Q_{1}\cap C^{c}_{\rho}})(0)
\leq N\int_{\rho}^{\infty}\frac{1}{r^{d+2-\gamma}}
\int_{0}^{r^{2}}\Big(\int_{|y|=r}f(s,y)\, \sigma_{r}(dy)\Big)\,dsdr,
$$
where $\sigma_{r}(dy)$ is the element of the surface area on $|y|=r$. Integrating by parts and using that   $\gamma<d+2$ 
(to ignore the integrated out terms) we get
$$
P_{\gamma}(fI_{Q_{1}\cap C^{c}_{\rho}})(0)\leq N
\int_{\rho}^{\infty}\frac{1}{r^{d+3-\gamma}}
\int_{\rho}^{r}\Big(\int_{0}^{\mu^{2}}
\Big(\int_{|y|=\mu}f(s,y)\,\sigma_{\mu}(dy)\Big)\,ds\Big)\,d\mu dr
$$
$$
\leq N
\int_{\rho}^{\infty}\frac{1}{r^{d+3-\gamma}}
\int_{0}^{r}\Big(\int_{0}^{r^{2}}
\Big(\int_{|y|=\mu}f(s,y)\,\sigma_{\mu}(dy)\Big)\,ds\Big)\,d\mu dr
$$
$$
=N\int_{\rho}^{\infty}\frac{1}{r^{d+3-\gamma}}I(r)\,dr ,
$$
where  
$$
I(r)=\int_{C_{r}}f(s,y)\,dyds.
$$
We use that  $I(r)\leq Nr^{ d+2-\beta }M_{\beta}f(0)$ and that $\gamma<\beta$. Then we see that 
\begin{equation}
                          \label{1.17.3}
P_{\gamma}(fI_{Q_{1}\cap C^{c}_{\rho}})(0)\leq N\rho^{ \gamma-\beta }M_{\beta}f(0).
\end{equation}

Next, by integrating by parts
we  obtain  that
$$
P_{\gamma}(fI_{Q_{2}\cap C^{c}_{\rho}})(0)
\leq
\int_{\rho^{2}}^{\infty}\frac{1}{s^{(d+2-\gamma)/2}}
\int_{|y|\leq\sqrt s}f(s,y)\,dyds
$$
$$
\leq N\int_{\rho^{2}}^{\infty}\frac{1}{s^{(d+4-\gamma)/2}}I(\sqrt s)\,ds=
N\int_{\rho}^{\infty}\frac{1}{r^{d+3-\gamma}}I(r)\,dr.
$$
This along with \eqref{1.17.3} proves  \eqref{1.17.2} and the lemma. \qed

To prove  Theorem \ref{theorem 6,6,1} we also need its homogeneous 
version for homogeneous Morrey
spaces $\dot E_{p,q,\beta} $.
Accordingly define
$$
\dot E^{1,2}_{p,q,\beta} =\{u:u,Du,D^{2}u,
\partial_{t}u\in  \dot E_{p,q,\beta} \},
$$
where $Du,D^{2}u,
\partial_{t}u$ are the Sobolev derivatives,
and provide $\dot E^{1,2}_{p,q,\beta} $ with an obvious norm.

\begin{lemma}
               \label{lemma 6,18.1}
Let   
$r\geq p$  and let
$$
2-\gamma< \beta \leq \frac{d}{p}+\frac{2}{q}
<2-\gamma+\frac{d}{r}.
$$ 
Then for any
$u\in \dot E^{1,2}_{p,q,\beta}$ its trace $u(0,\cdot)$ is uniquely defined
and  
\begin{equation}
                     \label{6.6.1}
\|D^{\gamma}u(0,\cdot)\|_{\dot E_{r,\beta+\gamma-2} }\leq N \|\partial_{t}u,
D^{2}u\|
_{ \dot E _{p,q,\beta}},
\end{equation} 
where the constant  $N$ depends only
on $d,p,q,r,\beta$.

\end{lemma}

Proof. Take $\zeta\in C^{\infty}_{0} (\bR^{d+1})$, such that $\zeta(0)=1\geq \zeta\geq0$, define $\zeta_{n}(t,x )=
\zeta(t/n^{2},x/n )$ and observe that,
as $n\to\infty$,
$$
 \big|\|
\partial_{t}(\zeta_{n}u)\|_{\dot E_{p,q,\beta}} -\|\zeta_{n}
\partial_{t}u\|_{\dot E_{p,q,\beta} }\big|\leq n^{-2}\sup |\partial_{t}\zeta|
\|
u\|_{\dot E_{p,q,\beta}} \to 0.
$$
Also
$$
 \big|\|
D(\zeta_{n}u)\|_{\dot E_{p,q,\beta}} -\|\zeta_{n}
Du\|_{\dot E_{p,q,\beta} }\big|\leq n^{-1}\sup |D\zeta|
\|
u\|_{\dot E_{p,q,\beta}} \to 0,
$$
$$
 \big|\|
D^{2}(\zeta_{n}u)\|_{\dot E_{p,q,\beta}} -\|\zeta_{n}
D^{2}u\|_{\dot E_{p,q,\beta} }\big|\leq n^{-2}\sup |D^{2}\zeta|
\|
u\|_{\dot E_{p,q,\beta}}
$$
$$
+2n^{-1}\sup |D\zeta|\|
Du\|_{\dot E_{p,q,\beta} } \to 0.
$$

It follows that it suffices to concentrate on $u$ that vanishes
for large $|t|+|x|^{2}$. In that case set $-f=\partial_{t}u+\Delta u$.
To further reduce our
problem observe that using translations show that it suffices
to prove that for any
  $\rho>0$,
$$
\rho^{\beta+\gamma-2}\dashnorm D^{\gamma}u(0,\cdot)\|_{L_{r}( B_{\rho})}\leq N\sup_{\rho_{1}\geq \rho}
\rho^{\beta }_{1}\dashnorm f\|
_{L_{p,q} (C_{\rho_{1}})   )}
$$
\begin{equation}
                     \label{6.19.1}
=  N\sup_{\rho_{1}\in[\rho,\rho+\rho_{2}]}
\rho^{\beta }_{1}\dashnorm f\|
_{L_{p,q} (C_{\rho_{1}})   },
\end{equation}
where $\rho_{2}$ is such that
$u(t,x)=0$ for $|t|+|x|^{2}\geq \rho_{2}^{2}$  
and the last equality is due to
$\beta \leq d/p+2/q$.

It is easy to pass to the limit
in \eqref{6.19.1} from smooth functions to arbitrary ones in
$W^{1,2}_{p,q}(C_{\rho_{2}})\supset
E^{1,2}_{p,q,\beta}(C_{\rho_{2}})$.
Therefore, we may assume that $u$ is smooth.
We thus reduced the general case to the task of proving the first estimate
in  \eqref{6.19.1} for smooth $u$ with compact support. One more reduction
is achieved by using the self-similarity
which shows that we only need to concentrate on $\rho=1$, that is
we only need to prove
\begin{equation}
                    \label{6,19,1}
 \dashnorm D^{\gamma}u(0,\cdot)\|_{L_{r}( B_{1})}\leq N\sup_{\rho \geq 1}
\rho^{\beta} \dashnorm f\|
_{L_{p,q} (C_{\rho })   )}
\end{equation}
for smooth $u$ with compact  support.

  Now   define
  $
g=|f|I_{C_{2}},h=|f|I_{C_{2}^{c}}$.
As it follows from the proof of Lemma~
\ref{lemma 6,17.1},
$$
|D^{\gamma}u(0,x)|\leq NG_{\gamma}(x)+NH_{\gamma}(x),
$$
where
$$
 (G_{\gamma},H_{\gamma})(x) =
\int_{0}^{\infty} 
\int_{\bR^{d}}P_{\gamma}(t,x-y)(g,h)(t,y)\,dydt.
$$

Estimates \eqref{1.17.2}
implies that
that for $|x|\leq 1$
$$
H_{\gamma}(x)\leq N\sup_{\rho>1}\rho^{\beta}
\dashint_{C_{\rho}(0,x)}h\,dydt
\leq N\sup_{\rho>1}\rho^{\beta}
\dashint_{ C_{2\rho}}h\,dydt
$$
$$
\leq N\sup_{\rho \geq1}
\rho^{\beta} \dashnorm f\|
_{L_{p,q}  (C_{\rho} )},
$$
where the last inequality is
due to H\"older's inequality.
Hence,
\begin{equation}
                     \label{6,8,3}
 \| H_{\gamma}\|_{L_{r}(  B_{1})}\leq N\sup_{\rho \geq1}
\rho^{\beta} \dashnorm f\|
_{L_{p,q}  (C_{\rho} )}. 
\end{equation}

Next, by Minkowski's inequality
$$
\|G_{\gamma}\|_{L_{r}(B_{1})}
\leq \int_{0}^{\infty}\Big(
\int_{B_{1}}\Big(\int_{\bR^{d}}
P_{\gamma}(t,x-y)g(t,y)\,dy\Big)^{r}\,dx\Big)^{1/r}
\,dt,
$$ 
where inside the integral with respect
to $t$ we have the norm of   convolution, so that by Young's
inequality this expression is dominated by
$$
\|g(t,\cdot)\|_{L_{p}}\|P_{\gamma}(t,\cdot)
\|_{L_{s}},
$$
where $1/s=1+1/r-1/p$ ($ \leq 1$ since $r\geq p$). We know that $\|P(t,\cdot)
\|_{L_{s}}=N(d)t^{\alpha}$ with
$\alpha=-\gamma/2+(d/2)(1/r-1/p)$, which after taking into account that 
$\|g(t,\cdot)\|_{L_{p}}=0$ for $t\geq4$, yields
$$
\|G\|_{L_{r}(B_{1})}
\leq N\int_{0}^{4}\|f(t,\cdot)\|_{L_{p}(B_{2})}t^{\alpha}\,dt.
$$
Now use H\"older's inequality along with the observation that
$\alpha q/(q-1)>-1$, due to 
the assumption that   $2-\gamma+d/r>
d/p+2/q$, to conclude that
$$
\|G\|_{L_{r}(B_{1})}
\leq N \|f \|_{L_{p,q}(C_{2})}.
$$
This and \eqref{6,8,3} prove
\eqref{6,19,1} and the lemma. \qed

{\bf Proof of Theorem \ref{theorem 6,6,1}}.  To prove \eqref{6,6.1}, it suffices to show that
for any $\rho\in (0, 1]$,
$\varepsilon>0$ 
\begin{equation}
                     \label{6,20.1}
I_{\rho}:=\rho^{\beta+\gamma-\mu}\dashnorm Du(0,\cdot)\|_{L_{r}(B_{\rho})}\leq N\varepsilon \|\partial_{t}u,
D^{2}u\|
_{ E _{p,q,\beta}}
+N\varepsilon^{-\mu/(2-\mu)} 
\|u\|
_{ E _{p,q,\beta}}.
\end{equation} 

By Corollary
\ref{corollary 6,19.1} with
$\epsilon=\varepsilon \rho^{ \gamma-\mu }$
in place of $\varepsilon$  we get
$$
 I_{\rho}\leq N\epsilon \|\partial_{t}u,
D^{2}u\|
_{ E _{p,q,\beta}}+N\big(\epsilon  \rho^{-2} +\epsilon^{-\kappa/(2-\kappa)}\rho^{(2\kappa-2\mu)/(2-\kappa)
}\big)\| u\|_{E_{p,q,\beta}}.
$$
For $\epsilon\leq \rho^{2-\mu}$ this yields (here we use that $ \kappa\leq \mu<2$)
$$
 I_{\rho}\leq N\epsilon \|\partial_{t}u,
D^{2}u\|
_{ E _{p,q,\beta}}+N\epsilon^{-\mu/(2-\mu)}\| u\|_{E_{p,q,\beta}}.
$$

In the remaining case $\rho^{2-\mu}<\epsilon$. In that case for $\zeta\in
C^{\infty}_{0}((-1,1)\times B_{2})$   such that $\zeta
=1$ on $C_{1}$ we have by Lemma \ref{lemma 6,18.1} that
\begin{equation}
                     \label{8,27.2}
I_{\rho}\leq \rho^{2-\mu} \rho^{\beta+\gamma-2}
\dashnorm D(\zeta u)(0,\cdot)\|_{L_{r}(B_{\rho})}\leq N\epsilon
\|\partial_{t}(\zeta u),D^{2}(\zeta u)
\|_{\dot E_{p,q,\beta}}.
\end{equation}
Owing to $\beta\leq d/p+2/q$,
 the last norm here is easily shown to be less than
$$
N\|\partial_{t} u,D^{2}u
\|_{  E_{p,q,\beta} }
+N\| u\|_{  E_{p,q,\beta}}.
$$ 
Therefore, \eqref{6,20.1} holds in this case as well and this proves
\eqref{6,6.1}. Estimate \eqref{8,27.1}
follows from \eqref{8,27.2} with
$\mu=2$ and $\epsilon=1$.
The theorem is proved. \qed

\end{document}